\documentclass[12pt, a4paper]{article}
\usepackage{a4wide}
\usepackage{setspace}
\usepackage[OT1]{fontenc}
\usepackage{amssymb,amsthm,amsmath}

\newcommand\X{\mathbf x}

\def\Im{{\rm\,Im} }

\newcommand\bY{\mathbf Y}
\newcommand\bR{{\mathbb R}}

\newcommand\y{{\mathbf y}}

\newcommand\wh{\widehat }
\newcommand\cov{\mathrm{Cov}}

\newcommand\e{{\mathbb E}}

\newcommand\var{{\rm Var}}
\newcommand\tr{{\mathrm{tr}}}

\newtheorem{theorem}{Theorem}[section]%
\newtheorem{corollary}[theorem]{Corollary}%
\newtheorem{remark}[theorem]{Remark}%

\begin{document}

\begin{center}
\Large Variance inequalities for quadratic forms with applications.  
\end{center}
\begin{center}
\large Pavel~Yaskov\footnote{Steklov Mathematical Institute of RAS, Russia\\
National University of Science and Technology MISIS, Russia\\
 e-mail: yaskov@mi.ras.ru\\The work is supported by the Russian Science Foundation
   under grant 14-21-00162.}
 \end{center}

\begin{abstract} We obtain variance inequalities for quadratic forms of weakly dependent random variables with bounded fourth moments. We also discuss two application. Namely, we  use these inequalities for deriving the limiting spectral distribution of a random matrix and estimating the long-run variance of a stationary time series.
\end{abstract}

\begin{center}
\end{center}

\section{Introduction}
Moment inequalities for quadratic forms  constitute a powerful tool in time series analysis and the random  matrix theory. In particular, they are used in the study of consistency and optimality properties of spectral density estimates (e.g., see Section V.4 in Hannan \cite{H}) as well as provide low-level conditions under which the limiting spectral distribution of a random matrix can be derived (e.g., see  Chapter 19 in  Pastur and Shcherbina \cite{PS}). 

When random variables $\{X_{i}\}_{i=1}^n$  are  independent, moment inequalities for quadratic forms $\sum_{i,j=1}^n a_{ij}X_iX_j$ are well-studied (e.g., see Lemma B.26 in Bai and Silverstein \cite{BS} and Chen \cite{C}). In the time series context, similar inequalities  were  obtained by many authors in connection with spectral density estimation and long-run variance estimation (e.g., see Chapter VI in White \cite{W}, Sections 6 and 7 in Wu and Xiao \cite{WX} and the references therein). In particular, high-order moment inequalities for causal time series were obtained by Wu and Xiao \cite{WX}. 

In the present paper we study variance inequalities for quadratic forms  $\sum_{i,j=1}^n a_{ij}X_iX_j$ of weakly dependent random variables $\{X_i\}_{i=1}^n$ with bounded fourth moments. 
Our  assumptions deal with covariances of $X_i$'s products up to the fourth order only and are closely related to the classical fourth-order cumulant condition  for a stationary time series (see Theorem V.4 in Hannan \cite{H} and Assumption A in Andrews \cite{A}). These assumptions can be easily verified under standard weak dependence conditions (e.g., strong mixing).  We also demonstrate how our  results can be applied in the random matrix theory and time series analysis.

The paper is structured as follows. Main results are given in Section 2. Section 3 is devoted to applications.
Section 4 contains all proofs.

\section{Main results}
Let $\{X_k\}_{k=1}^\infty $ be a sequence of centred random variables and let  $\{\varphi_k\}_{k=1}^\infty$  be a non-increasing sequence of non-negative numbers such that, for all $i<j<k<l$,
\begin{equation}\label{vphi}|\cov(X_i,X_jX_kX_l)|\leqslant \varphi_{j-i},\quad |\cov(X_iX_jX_k,X_l)|\leqslant \varphi_{l-k},
\end{equation}
\begin{equation}\label{vphi1}|\cov(X_iX_j,X_kX_l)|\leqslant \varphi_{k-j},\quad\text{and}\quad |\cov(X_i,X_{j})|\leqslant \varphi_{j-i}.
\end{equation}
Assumptions of this kind go back to R\'ev\'esz \cite{R} and studied by Koml\'os and R\'ev\'esz \cite{KR}, Gaposhkin \cite{G}, and  Longecker and Serfling \cite{LS} (see also M\'oricz \cite{M}, Section 4.3 and 4.4 in Doukhan et al. (2007)).

For simplicity, we will further assume that $\e X_k^2\leqslant 1$,  $k\geqslant 1$.
Define $\X_p=(X_1,\ldots,X_p)$ for $p\geqslant 1$, $\Phi_0=\sup\{\e X_1^4,\e X_2^4,\ldots \}$, and $ \Phi_1=\sum_{k=1}^\infty k\varphi_k.$ 
\begin{theorem}\label{t1}  
There is a  universal constant $C>0$ such that, for any $a\in\bR^p$ and all $p\times p$ matrices $A$ with zero diagonal, \[\e(\X_p^\top a)^4\leqslant C(\Phi_0+\Phi_1)(a^\top a)^2\;\;\text{and}\;\; \var(\X_p^\top A\X_p)\leqslant C(\Phi_0+\Phi_1 )\tr(A A^\top).\]
\end{theorem} 
A version of the first inequality in Theorem \ref{t1} is proved by  Koml\'os and R\'ev\'esz \cite{KR}, Gaposhkin \cite{G}, and  Longecker and Serfling \cite{LS}. The second inequality is new.

 Let now $\phi_k,$ $k\geqslant 1,$ satisfy 
\begin{equation}\label{phi}
\cov(X_i^2,X_j^2)\leqslant \phi_{j-i} \quad\text{for all }i<j.
\end{equation} 
Define $\Phi_2=\sum_{k=1}^\infty \phi_k$.
\begin{theorem}\label{t2}  There is a universal constant $C>0$ such that,  for all $p\times p$ matrices $A$,   
 \[\var(\X_p^\top A\X_p)\leqslant C(\Phi_0+\Phi_1+ \Phi_2)\tr(A A^\top).\]
\end{theorem} 

Let us give two  examples of $\{X_k\}_{k=1}^\infty$ that satisfy \eqref{vphi}, \eqref{vphi1}, and \eqref{phi}. \\
\noindent {\bf Example 1.} Let $\{X_k\}_{k=1}^\infty$ be a martingale difference sequence with  bounded 4th moments. Then \eqref{vphi} and \eqref{vphi1} hold for $\varphi_k=0,$ $k\geqslant 1,$ and $\Phi_1=0$. However, in general, there are no such $\phi_k,$ $k\geqslant 1$,  that \eqref{phi} holds and $\Phi_2<\infty$. This explains why we introduce two sets of coefficients $\{\varphi_k\}_{k=1}^\infty$ and $\{\phi_k\}_{k=1}^\infty$. If, in addition, $\{X_k^2-\e X_k^2\}_{k=1}^\infty$ is a martingale difference sequence, then, of course, $\phi_k=0$, $k\geqslant 1,$ and $\Phi_2=0.$

\noindent {\bf Example 2.} Let $\{X_k\}_{k=1}^\infty $ be  strongly mixing random variables  with mixing coefficients $(\alpha_k)_{k=1}^\infty$, zero mean,  and bounded moments of order $4\delta$ for some $\delta>1$. Then \eqref{vphi}--\eqref{phi} hold for $\varphi_k=\phi_k=C\alpha_k^{(\delta-1)/\delta}$ and large enough $C>0$ (e.g., see  Corollary A.2 in Hall and Heyde \cite{HH}). One can give similar bounds for other weak dependence conditions.

\begin{remark} We believe that higher order moment inequalities for quadratic forms $\X_p^\top A\X_p$ can be derived under similar conditions on $\cov(X_{i_1}\ldots X_{i_k},X_{i_{k+1}}\ldots X_{i_p})$ for $i_1< \ldots<i_p$, $k=1,\ldots,p-1$, and $p>4$.  However, the proofs are quite technical even in the case of the second-order inequalities and we leave this question for future research.
\end{remark}

Consider  the special case when $X_k,$ $k\geqslant 1$, are centred orthonormal random variables. In this case, \eqref{vphi} and \eqref{vphi1} reduce to 
\begin{equation*}|\e X_iX_jX_kX_l|\leqslant\min\{ \varphi_{j-i},\varphi_{k-j},\varphi_{l-k}\},\quad i<j<k<l.
\end{equation*}
Let $\y_p=(Y_1,\ldots,Y_p)$, where each $Y_j$  can be written as   $\sum_{k=1}^\infty a_{k}X_k$ in $L_2$ for some $a_{k}\in\bR,$ $k\geqslant 1,$  with $\sum_{k=1}^\infty a_k^2<\infty.$
\begin{corollary}\label{t4} Let $\Sigma_p=\e\y_p\y_p^\top$. Then there is $C>0$ such that, for any $p\times p$  matrix  $A$,  \[\var(\y_p^\top A\y_p)\leqslant C(\Phi_0+\Phi_1+\Phi_2)\tr(\Sigma_p A\Sigma_p A^\top).\] 
\end{corollary}  

If $\{X_k\}_{k=1}^\infty$ are independent standard normal variables and $A$ is a $p\times p$ symmetric matrix, then $\var(\y_p^\top A\y_p)=2\tr((\Sigma_p A)^2)$
(e.g., see Lemma 2.3 in Magnus \cite{JM}). Thus, Corollary \ref{t4} delivers an optimal estimate of the variance up to the factor $C(\Phi_0+\Phi_1+\Phi_2).$

\section{Applications}
In this section we discuss two applications of the obtained inequalities.

 Our first application will be in  the random matrix theory.
For each $p,n\geqslant 1$, let $\bY_{pn}$ be a $p\times n$ random matrix whose columns are independent copies
of $\y_p$, where $\y_p$ is given either in Corollary \ref{t4}, or $\y_p=\X_p$ for $\X_p$ from Theorem \ref{t2}. 
\begin{theorem}\label{rm}
Let  $\Phi_0,\Phi_1,\Phi_2<\infty$. If the following conditions hold
\begin{itemize}
\item[$(1)$] $p=p(n)$ is such that $p/n\to c$ for some $c>0,$
\item[$(2)$] the spectral norm of $\Sigma_{p}=\e \y_p\y_p^\top $ is bounded over $p$,
\item[$(3)$]  the empirical spectral distribution of $\Sigma_p$'s eigenvalues  has a  weak  limit $P(d\lambda)$,
\end{itemize} 
then, with probability one, the empirical spectral distribution of $n^{-1}\bY_{pn}\bY_{pn}^\top$'s eigenvalues has a  weak  limit whose Stieltjes transform $m = m(z)$ satisfies
\[m(z)=\int_0^\infty\frac{P(d\lambda)}{\lambda(1-c-czm(z))-z},\quad z\in\mathbb C,\quad \Im(z)>0.\]
\end{theorem}

The next application concerns long-run variance estimation. First, let us recall the generic form of the central limit theorem for  a weakly stationary time series $(X_t)_{t=-\infty}^\infty$:
\[\sqrt{n}(\overline{X}_n-\mu)\stackrel{d}\to\mathcal N(0,\sigma^2),\quad n\to\infty,\]
where $\overline{X}_n=n^{-1}\sum_{t=1}^n X_t$, $\e X_t=\mu$, and $\sigma^2$ is the long-run variance of $(X_t)_{t=-\infty}^\infty$, i.e.
\[\sigma^2=\sum_{j=-\infty}^\infty\cov(X_t,X_{t+j}).\]
This theorem can be proved under different weak dependence assumptions (e.g., see the books of Doukhan et al. \cite{D} and Bulinski and Shashkin \cite{BS}). In statistical applications, this theorem takes the form 
\[\frac{\sqrt{n}(\overline{X}_n-\mu)}{\wh\sigma_n}\stackrel{d}\to\mathcal N(0,1),\quad n\to\infty,\]
where $\wh\sigma^2_n$ is a consistent estimator of $\sigma^2.$ Recall also that $\sigma^2$ can be written as  $\sigma^2=2\pi f(0),$ where $f=f(x),$ $x\in[-\pi,\pi),$  is the spectral density of $(X_t)_{t=-\infty}^\infty$. Therefore, long-run variance estimation is closely related to spectral density estimation.

A number of papers is devoted to the study of consistency and optimality properties of long-run variance estimators (e.g., see Andrews \cite{A}, Hansen \cite{H1}, de Jong and Davidson \cite{JD}, and Jansson \cite{J} among others). When $\mu=0$, a typical estimator has the form
\begin{equation}\label{sig}
\wh\sigma^2_n=\frac{1}{n}\sum_{s,t=1}^n K\bigg(\frac{|s-t|}{m}\bigg)X_sX_t,
\end{equation} 
where $K=K(x),$ $x\geqslant 0,$ is a kernel function. Standard assumptions on $K=K(x)$ include

$(a)$ $K(0)=1$, $K$ is  continuous at $x=0$, and $\sup_{x\geqslant 0}|K(x)|<\infty$,

$(b)$ $\int_{0}^\infty \bar{K}^2(x)\,dx<\infty$ for $\bar{K}(x)=\sup_{y\geqslant x}|K(y)|,$

$(c)$ $k_q=\lim_{x\to 0+}x^{-q}(K(x)-1)$ exists for some $q> 0$.

Assumptions $(a)$--$(b)$ are inspired by Assumption 3 of Jansson \cite{J}. However, $(b)$ is weaker than Assumption 3(ii) in \cite{J}, where  the integrability of $\bar{K}$ is assumed.
 To our knowledge, the weakest alternative to $(b)$ considered in the literature
is the integrability of  $K^2$. However, as discussed  in Jansson \cite{J}, many previous results (Andrews \cite{A}, Hansen \cite{H1}, among others)  are incorrect as they stated and need stronger alternatives to  the integrability of  $K^2$.  Assumption $(c)$ is classical and goes back to Parzen \cite{Pa} (see also Andrews \cite{A}). 

Let further $(X_t)_{t=-\infty}^\infty$ be a centred weakly stationary time series that satisfies conditions from Section 2 (in particular, $\e X_t^2\leqslant 1$).  Our first result is the consistency of $\wh\sigma^2_n$. 

\begin{theorem}\label{lv} Let $K=K(x)$ satisfy $(a)$--$(b)$. If $\Phi_0,\Phi_1,\Phi_2$ are finite, then 
$\wh\sigma^2_n\to \sigma^2$ in mean square as $m,n\to\infty$ and $m=o(n)$. 
\end{theorem}

The dependence $m=o(n)$ is optimal. It can be  seen by taking a Gaussian white noise  $(X_t)_{t=-\infty}^\infty$  and showing that $\var(\wh\sigma^2)\not\to 0$ when $m/n\not \to0$ due to the variance formula for Gaussian quadratic forms given in the end of Section 2. Andrews \cite{A} following Hannan \cite{H} proved consistency of $\wh\sigma^2_n$ under the cumulant condition 
\begin{equation}\label{cum}
\sum_{j,k,l=1}^\infty \sup_{t\geqslant 1}|\kappa(X_t,X_{t+j},X_{t+k},X_{t+l})|<\infty.
\end{equation}
Here $\kappa(X_i,X_{j},X_{k},X_{l})$ is the fourth-order cumulant that is equal to
 \[\e X_iX_{j}X_{k}X_{l}-\e X_iX_{j} \e X_{k}X_{l} -
\e X_iX_{k}\e X_{j}X_{l}- \e X_iX_{l}\e X_{k}\e X_{j}\]
when each $X_t$ has  zero mean. By Lemma 1 of Andrews \cite{A}, \eqref{cum} holds when  $(X_t)_{t=-\infty}^\infty$  is a strongly mixing sequence with mixing coefficients $(\alpha_k)_{k=1}^\infty$ satisfying \[\sum_{k=1}^\infty k^2 \alpha_k^{(\delta-1)/\delta}<\infty\] and bounded moments of order $4\delta$ for some $\delta>1$. By Example 2, our Theorem \ref{lv} is applicable whenever $\sum_{k=1}^\infty k \alpha_k^{(\delta-1)/\delta}<\infty$.

The cumulant condition allows to calculate the limit of the mean squared error (MSE) of $\wh\sigma_n^2$ explicitly. We can not do it under our assumptions. However, we can give an upper bound for MSE which is very similar to the exact limit (see Proposition 1 in Andrews \cite{A}).

\begin{theorem}\label{lv1} Under conditions of Theorem \ref{lv}, let $(c)$ hold for some $q>0$ and
 \[\text{the series }\Gamma_q=\sum_{j=1}^\infty j^q\cov(X_t,X_{t+j})\text{ converges absolutely.}\]Then there is an absolute constant $C>0$ such that,  as $m,n\to\infty,$
\[\e|\wh\sigma^2_n-\sigma^2|^2\leqslant C(\Phi_0+\Phi_1+\Phi_2)\frac{m}{n}\int_{0}^\infty \bar{K}^2(x)\,dx+\frac{4(k_q\Gamma_q)^2}{m^{2q}}+o(m^{-2q})+O(n^{-1}).\]
\end{theorem}

\section{Proofs}
Below we assume that $\Phi_0,\Phi_1,\Phi_2$ are finite  otherwise all bounds become trivial.

\noindent {\bf Proof of Theorem \ref{t1}.} To prove the first inequality, we reproduce the proof given in  Gaposhkin \cite{G} with the only difference that we derive explicit constants in his inequality. 
Note first that, as $\e X_i=0$ for all $i\geqslant 1$, it follows from \eqref{vphi} that 
\begin{equation}
\label{vphi2}
|\e X_iX_jX_kX_l|\leqslant\min\{\varphi_{j-i},\varphi_{l-k}\},\quad i<j<k<l.
\end{equation}

 Write $a=(a_1,\ldots,a_p)$. Using  Lemma 1 in  Moricz \cite{M} with $p=2$ and $r=4$, we get  \[|(\X_p^\top a)^4-24 T|\leqslant C_0(S^4+|\X_p^\top a|^3S),\]
where $C_0>0$ is a universal constant,
\[T=\sum_{i<j<k<l}a_ia_ja_ka_lX_iX_jX_kX_l,\quad S=\Big(\sum_{i=1}^pa_i^2X_i^2\Big)^{1/2},\]
hereinafter $i,j,k,l$ are any numbers in $\{1,\ldots,p\}$.
By H\"older's inequality,
\begin{align*}
\e(\X_p^\top a)^4\leqslant 24 \e T+C_0\big(\e S^4 +(\e S^4)^{1/4}(\e |\X_p^\top a|^4)^{3/4}).
\end{align*}
 By  \eqref{vphi2},
\[
|\e T|\leqslant\sum_{ i<j<k<l}|a_ia_ja_ka_l| \min\{\varphi_{j-i},\varphi_{l-k}\}\leqslant \frac{1}{4}\sum_{ i<j<k<l}(a_i^2+a_j^2)(a_k^2+a_l^2) \min\{\varphi_{j-i},\varphi_{l-k}\}. \]
We estimate only the term 
\[J=\sum_{ i<j<k<l}a_i^2a_k^2\min\{\varphi_{j-i},\varphi_{l-k}\}. \]
Other terms with $a_j^2a_k^2$, $a_i^2a_l^2$, and $a_j^2a_l^2$ instead of $a_i^2a_k^2$ can be estimated similarly. We have
\[J\leqslant \sum_{ i<k}a_i^2a_k^2\sum_{q,r=1}^\infty \min\{\varphi_{q},\varphi_{r}\}\]
and 
\begin{equation}
\label{qr}
\sum_{q,r=1}^\infty \min\{\varphi_{q},\varphi_{r}\}\leqslant
\sum_{q=1}^\infty \Big( q\varphi_q+\sum_{r=q+1}^\infty  \varphi_r \Big)=\Phi_1+\sum_{r=2}^\infty\sum_{q=1}^{r-1}\varphi_r\leqslant 2\Phi_1.
\end{equation}
As a result,  \[J\leqslant 2\Phi_1\sum_{ i<k}a_i^2a_k^2\leqslant \Phi_1\Big(\sum_{ i=1}^pa_i^2\Big)^2\quad\text{and}\quad |\e T|\leqslant \Phi_1\|a\|^4.\]

Let us also note that 
\[\e S^4=\sum_{i,j=1}^pa_i^2 a_j^2\e X_i^2X_j^2\leqslant \Phi_0\|a\|^4.\]
Combining the above estimates, we infer that
\[\e(\X_p^\top a)^4\leqslant (24+C_0) (\Phi_0+\Phi_1)\|a\|^4 +C_0[(\Phi_0+\Phi_1)\|a\|^4]^{1/4}[\e(\X_p^\top a)^4]^{3/4}.\]
Put $R=[\e(\X_p^\top a)^4/(\Phi_0+\Phi_1)]^{1/4}/\|a\|$. Then $R^4\leqslant 24+C_0+C_0R^3.$ Therefore, $R\leqslant R_0$, where $R_0>0$ is the largest positive root of the equation $x^4=24+C_0+C_0x^3.$ Finally, we obtain
\[
\e|\X_p^\top a|^4\leqslant R_0^4(\Phi_0+\Phi_1)\|a\|^4.
\]

We now verify the second inequality.  First, note that, for $i<j,$ 
\[|\cov( X_i,X_j)|=|\e X_iX_j|\leqslant \sqrt{\e X_i^2\e X_j^2}\leqslant 1\]
In addition, for $i<j<k<l,$ \begin{equation}\label{fi}
|\cov(X_iX_j,X_kX_l)|\leqslant 2\min\{\varphi_{j-i},\varphi_{k-j},\varphi_{l-k}\}\quad\text{and}\quad I\leqslant 2\min\{\varphi_{j-i},\varphi_{l-k}\},
\end{equation}  
where $I$ is equal to $|\cov(X_iX_k,X_jX_l)|$ or  $|\cov(X_iX_l,X_jX_k)|$. Indeed, by \eqref{vphi1} and  \eqref{vphi2}, 
\begin{align*}|\cov(X_iX_j,X_kX_l)|&\leqslant\min\{\varphi_{k-j},|\e X_iX_jX_kX_l|+| \e X_iX_j\e X_kX_l|\}\\
&\leqslant \min\{\varphi_{k-j},2\min\{\varphi_{j-i},\varphi_{l-k}\}\}\\
&\leqslant 2\min\{\varphi_{j-i},\varphi_{k-j},\varphi_{l-k}\},
\end{align*}
and, by  the monotonicity of $\varphi_k$, 
\begin{align*}
|\cov(X_iX_k,X_jX_l)|&\leqslant |\e X_iX_k\e X_jX_l|+|\e X_iX_j X_kX_l |\\
&\leqslant\min\{\varphi_{k-i},\varphi_{l-j}\}+\min\{\varphi_{j-i},\varphi_{l-k}\}\\
&\leqslant 2\min\{\varphi_{j-i},\varphi_{l-k}\}.
\end{align*}
A similar bound holds for $\cov(X_iX_l,X_jX_k).$

Let $A=(a_{ij})_{i,j=1}^p$ and $a_{ii}=0$, $1\leqslant i\leqslant p$. Set $B=(A^\top+A)/2$. Then $\X_p^\top A\X_p=\X_p^\top B\X_p$ and
\begin{equation}\label{BB}
\tr(BB^\top)=\sum_{i,j=1}^p \Big(\frac{a_{ij}+a_{ji}}2\Big)^2\leqslant
\sum_{i,j=1}^p \frac{a_{ij}^2+a_{ji}^2}2=\sum_{i,j=1}^p a_{ij}^2=\tr(AA^\top).
\end{equation}
Hence, we may assume w.l.o.g. that $A=A^\top.$  Then 
\begin{align*}
\var(&\X_p^\top A\X_p)=4
\var\Big(\sum_{i=1}^{p-1} X_i \sum_{k=i+1}^p a_{ik} X_k\Big)=4\sum_{i=1}^{p-1} \var\Big( X_i\sum_{k=i+1}^p a_{ik} X_k\Big)+\\
&+8\sum_{i<j}\cov\Big( X_i\sum_{k=i+1}^p a_{ik} X_k,X_j\sum_{k=j+1}^p a_{jk} X_k\Big)=4I_1+8I_2+8I_3+8I_4,
\end{align*}
where
\begin{align*}
I_1&=\sum_{i=1}^{p-1} \var\Big( X_i\sum_{k=i+1}^p a_{ik} X_k\Big),\\
I_2&=\sum_{i<j}\cov\Big( X_i \sum_{k=i+1}^{j-1} a_{ik} X_k,  X_j \sum_{k=j+1}^p a_{jk} X_k\Big),\\
I_3&=\sum_{i<j}a_{ij}\cov\Big( X_i X_j,X_j \sum_{k=j+1}^p a_{jk} X_k\Big),\\
I_4&=\sum_{i<j}\cov\Big( X_i\sum_{k=j+1}^p a_{ik} X_k,X_j\sum_{k=j+1}^p a_{jk} X_k\Big),
\end{align*}
and the sums over the empty set are zeros.
\\
{\it Control of $I_1$.} By the Cauchy-Schwartz inequality and the first inequality in  Theorem \ref{t1},
\begin{align*}
I_1\leqslant &
\sum_{i=1}^{p-1} \sqrt{\e X_i^4}\Big(\e\Big|\sum_{k=i+1}^p a_{ik} X_k\Big|^4\Big)^{1/2}\\
\leqslant&
 C(\Phi_0+\Phi_1) \sum_{i=1}^{p-1} \sum_{k=i+1}^p a_{ik}^2=C(\Phi_0+\Phi_1)\frac{\tr(A^2)}2
\end{align*}
\\
{\it Control of $I_2$.} By the Cauchy inequality and \eqref{fi},
\begin{align*} I_2\leqslant&\sum_{i<k<j<l}|a_{ik}a_{jl}|\,|\cov(X_iX_k,X_jX_l)|\\
\leqslant& 2\sum_{i<k<j<l}|a_{ik}a_{jl}|\min\{\varphi_{k-i},\varphi_{j-k},\varphi_{l-j}\}\\
\leqslant & I_5+I_6 ,\end{align*}
where
\[ I_5=\sum_{i<k<j<l}
a_{ik}^2\min\{\varphi_{j-k},\varphi_{l-j}\},\quad 
I_6=\sum_{i<k<j<l}
a_{jl}^2\min\{\varphi_{k-i},\varphi_{j-k}\}.
\]
Additionally, by \eqref{qr},
\begin{align*}
I_5\leqslant&\sum_{i<k}a_{ik}^2\sum_{q,r=1}^\infty \min\{\varphi_q,\varphi_r\}\leqslant   \frac{\tr(A^2)}2\,(2\Phi_1)= \tr(A^2)\Phi_1.
\end{align*}
We similarly derive that $I_6\leqslant\tr(A^2)\Phi_1.$  Hence, $I_2\leqslant 2\tr(A^2)\Phi_1 $.
\\
{\it Control of $I_3$.} By the Cauchy-Schwartz inequality and the first inequality in  Theorem \ref{t1},
\begin{align*}
I_3=&\sum_{j=2}^{p-1}\cov\Big(X_j\sum_{i=1}^{j-1} a_{ij} X_i,X_j\sum_{k=j+1}^p a_{jk} X_k\Big)\\
&\leqslant \sum_{j=2}^{p-1}\Big(\e  X_j^2\Big|\sum_{i=1}^{j-1} a_{ij} X_i\Big|^{2}\Big)^{1/2}\Big(\e X_j^2\Big|\sum_{k=j+1}^p a_{jk} X_k\Big|^2\Big)^{1/2}\\
&\leqslant\sum_{j=2}^{p-1}\sqrt{\e  X_j^4}\bigg[\e\Big(\sum_{i=1}^{j-1} a_{ij} X_i\Big)^4\e\Big(\sum_{k=j+1}^p a_{jk} X_k\Big)^4\bigg]^{1/4}\\&\leqslant\sqrt{C(\Phi_0+\Phi_1)}(I_7+I_8)/2,
\end{align*}
where
\[I_7=\sum_{j=2}^{p-1}\bigg[\e\Big(\sum_{i=1}^{j-1} a_{ij} X_i\Big)^4\bigg]^{1/2},\quad I_8=\sum_{j=2}^{p-1}\bigg[\e\Big(\sum_{k=j+1}^p a_{jk} X_k\Big)^4\bigg]^{1/2}.\]
By the first inequality in Theorem \ref{t1},
\[I_7\leqslant K\sum_{j=2}^{p-1}\sum_{i=1}^{j-1} a_{ij}^2\leqslant\frac{K\tr(A^2)}{2},\quad 
I_8\leqslant K\sum_{j=2}^{p-1}\sum_{k=j+1}^p a_{jk}^2\leqslant\frac{K\tr(A^2)}{2},\]
where $K=\sqrt{C(\Phi_0+\Phi_1)}.$ As a result, $I_3\leqslant C(\Phi_0+\Phi_1)\tr(A^2)/2$.
\\
{\it Control of $I_4$.} We have $I_4=I_9+I_{10}+I_{11}$, where
\[I_9=\sum_{i<j<k}\cov(a_{ik }X_i X_k,a_{jk} X_jX_k),\;\; I_{10}=
\sum_{i<j<k<l}a_{ik}a_{jl}\cov( X_iX_k,X_{j}X_l),\]\[
I_{11}=\sum_{i<j<k<l}a_{il}a_{jk}\cov( X_iX_l,X_{j}X_k).\]
By the first inequality in Theorem \ref{t1},
\[I_9=\frac{1}{2}\sum_{k=3}^p \var\Big(X_k\sum_{i=1}^{k-1}a_{ik}X_i\Big)-
\frac{1}{2}\sum_{k=3}^p \sum_{i=1}^{k-1}\var( a_{ik}X_iX_k) \leqslant \]
\[\leqslant\frac{1}{2} \sum_{k=3}^p\bigg[ \e X_k^4\e\Big(\sum_{i=1}^{k-1}a_{ik}X_i\Big)^4\bigg]^{1/2}\leqslant C(\Phi_0+\Phi_1)
\sum_{k=3}^p \sum_{i=1}^{k-1}\frac{a_{ik}^2}2\leqslant C(\Phi_0+\Phi_1)\,\frac{\tr(A^2)}{4}.
\]

Let us now estimate  $I_{10}$ and $I_{11}$. By \eqref{fi},
\begin{align*}
I_{10}\leqslant 2\sum_{i<j<k<l}|a_{ik}a_{jl}| \min\{\varphi_{j-i},\varphi_{l-k}\}\quad\text{and}\quad
I_{11}\leqslant 2\sum_{i<j<k<l}|a_{il}a_{jk}| \min\{\varphi_{j-i},\varphi_{l-k}\}.
\end{align*}  
By the Cauchy inequality,  $I_{10}\leqslant I_{12}+I_{13}$ and  $I_{11}\leqslant I_{14}+I_{15}$ with
\[I_{12}=\sum_{i<j<k<l}a_{ik}^2 \min\{\varphi_{j-i},\varphi_{l-k}\},\quad I_{13}=\sum_{i<j<k<l}a_{jl}^2 \min\{\varphi_{j-i},\varphi_{l-k}\},\]
\[I_{14}=\sum_{i<j<k<l}a_{il}^2 \min\{\varphi_{j-i},\varphi_{l-k}\},\quad I_{15}=\sum_{i<j<k<l}a_{jk}^2 \min\{\varphi_{j-i},\varphi_{l-k}\},\]
As previously, we have 
\begin{align*}
I_{12}&\leqslant\sum_{i<k}a_{ik}^2\sum_{q,r=1}^\infty \min\{\varphi_q,\varphi_r\}\leqslant \tr(A^2)\Phi_1.
\end{align*}
By the same arguments, $I_{13}$, $I_{14},$ and $I_{15}$ can be bounded from above by $\tr(A^2)\Phi_1$. 
Thus, we conclude that $I_{10}+I_{11}\leqslant  4\tr(A^2)\Phi_1$.

Combining all above estimates, we get
$\var(\X_p^\top A\X_p)\leqslant C(\Phi_0+\Phi_1)\tr(A^2)$ for a universal constant $C>0.$
Q.e.d.
\\
\noindent{\bf Proof of Theorem \ref{t2}.}
Let $A=(a_{ij})_{i,j=1}^p$ and $D$ be the $p\times p$ diagonal matrix with diagonal entries $a_{11},\ldots,a_{pp}$.
By Theorem \ref{t1},
 \[\var(\X_p^\top (A-D)\X_p) \leqslant C(\Phi_0+\Phi_1)\tr((A-D)(A-D)^\top).\] 
In addition,  $\var(\X_p^\top A\X_p)\leqslant 2\var(\X_p^\top D\X_p)+2\var(\X_p^\top (A-D)\X_p) .$ 
Noting that \[\tr (AA^\top)=\tr((A-D)(A-D)^\top)+\tr(D^2),\] we only need to bound $\var(\X_p^\top D\X_p)$ from above by $\tr(D^2)$ up to a constant factor. Write $D=D_1-D_2$, where $D_i$ are diagonal matrices with non-negative diagonal entries and $\tr(D^2)=\tr(D_1^2)+\tr(D_2^2)$.  By the Cauchy inequality, $\var(\X_p^\top D\X_p)\leqslant 2\sum_{i=1}^2\var(\X_p^\top D_i\X_p).$ Hence,   
 we may assume w.l.o.g. that diagonal elements of $D$ are non-negative.

We see that 
\begin{align*}
\var(\X_p^\top D\X_p )=\var\Big(\sum_{i=1}^pa_{ii}X_i^2 \Big)&=\sum_{i=1}^p a_{ii}^2\var( X_i^2)+\sum_{i\neq j}a_{ii}a_{jj}\cov( X_i^2, X_j^2)\\
&\leqslant \Phi_0\sum_{i=1}^n a_{ii}^2+\sum_{i\neq j}a_{ii}a_{jj}\phi_{|i-j|}
\end{align*}
and, as a result,
\begin{align*}
\var(\X_p^\top D\X_p )\leqslant&\Phi_0\tr(D^2)+\sum_{i\neq j}\frac{a_{ii}^2+a_{jj}^2}2 \phi_{|i-j|}\\
\leqslant&\Phi_0\tr(D^2)+\sum_{i=1}^p a_{ii}^2\sum_{j: j\neq i}\phi_{|i-j|}\\
\leqslant& 2\tr(D^2)\Big(\Phi_0+\sum_{k=1}^\infty \phi_k\Big)= 2(\Phi_0+\Phi_2)\tr (D^2).
\end{align*}
Combining the above bounds, we get the desired inequality. Q.e.d.\\
{\bf Proof of Corollary \ref{t4}.}
By the definition of $\y_p$, $\Gamma_n \X_n\to \y_p $ in probability and in mean square as $n\to\infty$ for some $p\times n$ matrices $\Gamma_n$ and $\X_n=(X_1,\ldots,X_n)$. Since $X_k, k\geqslant 1,$ are orthonormal, we have  $\Gamma_n\Gamma_n^\top=\e (\Gamma_n \X_n)(\Gamma_n \X_n)^\top\to \e \y_p\y_p^\top=\Sigma_p$, \[\X_n^\top(\Gamma_n^\top A \Gamma_n)\X_n=(\Gamma_n \X_n)^\top A( \Gamma_n \X_n)\to \y_p^\top A\y_p\quad\text{in probability},\]
and \[\e \X_n^\top( \Gamma_n^\top A\Gamma_n) \X_n=\tr(\Gamma_n^\top A \Gamma_n)=\tr( \Gamma_n\Gamma_n^\top A)\to \tr(\Sigma_pA)=\e \y_p^\top A\y_p\] as  $n\to\infty.$   We need the following version of Fatou's lemma:
\begin{center}
 if $\xi_n\to \xi$ in probability, then $\e |\xi|\leqslant \varliminf\limits_{n\to\infty}\e |\xi_n|.$
\end{center}By this lemma and  Theorem \ref{t2},
\begin{align*}
\e|\y_p^\top A\y_p-\tr(\Sigma_p A)|^2\leqslant & \varliminf\limits_{n\to\infty} \e|\X_n^\top( \Gamma_n^\top A\Gamma_n) \X_n-\tr(\Gamma_n^\top A\Gamma_n)|^2\\
\leqslant&\varliminf\limits_{n\to\infty} C(\Phi_0+\Phi_1+\Phi_2) \tr(\Gamma_n^\top A\Gamma_n\Gamma_n^\top A^\top \Gamma_n)
\end{align*}
Note that $\tr(\Gamma_n^\top A\Gamma_n\Gamma_n^\top A^\top \Gamma_n)=\tr(\Gamma_n\Gamma_n^\top A\Gamma_n\Gamma_n^\top A^\top)\to\tr(\Sigma_pA\Sigma_pA^\top ).$ Q.e.d.

\noindent{\bf Proof of Theorem \ref{rm}.} Denote  the spectral norm of a  matrix $A$ by  $\|A\|$. Recall that $\|A\|=\sqrt{\|AA^\top\|}=\sqrt{\|A^\top A\|}$. In addition, let $A^{1/2}$ be  the principal square root of a square positive semi-definite matrix $A$.  By Theorem 1.1 in Bai and Zhou \cite{BZ}, we will prove the theorem by checking that $\var(\y_p^\top A_p\y_p)=o(p^2)$ as $p\to\infty$ for any sequence $(A_p)_{p=1}^\infty$ with $\|A_p\|=O(1)$, where $A_p$ is a $p\times p$ matrix.

First, let $\y_p$ be as in Corollary \ref{t4}. Then \[\tr(\Sigma_pA_p\Sigma_pA_p^\top)=\tr(\Sigma_p^{1/2}A_p\Sigma_pA_p^\top \Sigma_p^{1/2})=\tr(Q \Sigma_pQ^\top)\]
with  $Q=\Sigma_p^{1/2}A_p$. If $I_p$ is the $p\times p$ identity matrix, then $\|\Sigma_p\|I_p-\Sigma_p$ is positive semi-definite and, as a result, $Q(\|\Sigma_p\|I_p-\Sigma_p) Q^\top$ is positive semi-definite for any $Q$. Hence,
\[\tr(Q^\top \Sigma_pQ)\leqslant\|\Sigma_p\|\tr(QQ^\top )=\|\Sigma_p\|\tr(Q^\top Q)=\|\Sigma_p\|\tr(A_p^\top \Sigma_p A_p)\leqslant \|\Sigma_p\|^2\tr(A_p^\top A_p)
\]
and $\tr(A_p^\top A_p)\leqslant \|A_p\|^2p$.  Therefore, by Corollary \ref{t4}, \[\var(\y_p^\top A_p\y_p)\leqslant C(\Phi_0+\Phi_1+\Phi_2)\|A_p\|^2\|\Sigma_p\|^2p=o(p^2)\]
whenever $\|A_p\|=O(1).$ The case $\y_p=\X_p$ with $\X_p$ as in Theorem \ref{t2} can be considered along the same lines due to the inequality  $\tr(A_pA_p^\top)\leqslant \|A_p\|^2p$. Q.e.d.

\noindent {\bf Proof of Theorem \ref{lv}.} Since $\Phi_1<\infty$, we have \begin{equation}
\label{covc} \sum_{j=1}^\infty j|C(j)|<\infty\end{equation} and $\sigma^2$ is well-defined, where $C(j)=\cov(X_t,X_{t+j}),$ $j\in\mathbb Z$. We also have the bias-variance decomposition $\e (\wh\sigma^2-\sigma^2)^2=\var(\wh\sigma^2)+(\e\wh\sigma^2-\sigma^2)^2.$

First, let us estimate the bias term $\e\wh\sigma^2-\sigma^2.$  Using  $K(0)=1,$ we get
\begin{align*}
\e\wh\sigma^2-\sigma^2=&\sum_{j=-n}^n \bigg(1-\frac{|j|}{n}\bigg)K\Big(\frac{|j|}m\Big)C(j)-\sum_{j=-\infty}^\infty C(j)\\
=&2\sum_{j=1}^n \bigg(K\Big(\frac{j}m\Big)-1\bigg)C(j)-\frac{2}n\sum_{j=1}^nK\Big(\frac{j}m\Big)jC(j)-2\sum_{j>n} C(j)
\end{align*}
Now, setting $M=\sup_{x\geqslant 0}|K(x)|,$
\[\Big|\sum_{j=1}^n K(j/m)jC(j)\Big|\leqslant M\sum_{j=1}^\infty  j|C(j)|=O(1).\]
Additionally,
\[\Big|\sum_{j>n} C(j)\Big|\leqslant \frac{1}{n}\sum_{j>n} j|C(j)|=o(1/n). \]
Combining these relations yields 
\begin{equation}
\label{bias}\e\wh\sigma^2-\sigma^2=2\sum_{j=1}^n (K(j/m)-1)C(j) +O(1/n)=o(1),
\end{equation}
 where the last equality follows from $(a)$ and \eqref{covc}.

Now, consider the variance term $\var(\wh\sigma^2)$. By Theorem \ref{t2}, 
\[\var(\wh\sigma^2)\leqslant \frac{C_0}{n^2}\sum_{s,t=1}^n K^2\bigg(\frac{|s-t|}m\bigg), \]
where $C_0=C(\Phi_0+\Phi_1+\Phi_2)$ with $C$ given in Theorem \ref{t2}. Using that  $\bar K(x)=\sup_{y\geqslant x}|K(y)|$ is a non-decreasing function in $L_2(\bR)$, we derive 
\[\frac{1}{mn}\sum_{s,t=1}^n K^2\bigg(\frac{|s-t|}m\bigg)\leqslant\frac{1}{m}+ \frac{2}{m}\sum_{j=1}^n \bar K^2\Big(\frac{j}m\Big)\leqslant \frac{1}{m}+2\int_0^\infty \bar K^2(x)\,dx.\]
As a result,
\begin{equation}\label{var}
\var(\wh\sigma^2)\leqslant \frac{C_0}{n}+\frac{2C_0m}{n}\int_0^\infty \bar K^2(x)\,dx
\end{equation}
and   $\var(\wh\sigma^2)=o(1)$ whenever $m,n\to\infty$ and $m/n\to0.$

Combining the above bounds for the bias and the variance, we finish the proof. Q.e.d.

\noindent {\bf Proof of Theorem \ref{lv1}.} The proof follows the same line as the proof of Theorem \ref{lv}. We only need to note the following. If $(c)$ holds for some $q>0$. Then, by $(c)$ and  the boundedness of $K$, $x^{-q}(K(x)-1)$ is bounded on $\bR_+$. Therefore, by $(a)$, $(c)$, and the absolute convergence of $\sum_{j\geqslant 1}j^qC(j)$, 
 \[m^q\sum_{j=1}^n (K(j/m)-1)C(j)=\sum_{j=1}^n \frac{K(j/m)-1}{(j/m)^q}\,j^qC(j)=k_q\sum_{j=1}^\infty j^qC(j)+o(1).\] 
By  \eqref{bias} and  \eqref{var}, the latter yields the desired bound. Q.e.d.

\end{document}